\documentclass[12pt]{article}
\usepackage{amsmath}
\usepackage{amssymb}
\usepackage{tabularx}
\usepackage{enumerate}
\usepackage{graphicx}
\usepackage{color}

\newtheorem{thm}{Theorem}[section]
\newtheorem{lem}[thm]{Lemma}

\newtheorem{rmk}{Remark}[section]
\newtheorem{defi}{Definition}[section]

\newtheorem{pppp}{Proof}

\newcommand{\qed}{\hspace{1em}\mbox{\raisebox{0.65ex}{\fbox{}}}}

\numberwithin{equation}{section}

\newcommand{\be}{\begin{equation}}
\newcommand{\ee}{\end{equation}}
\newcommand\bes{\begin{eqnarray}} \newcommand\ees{\end{eqnarray}}
\newcommand{\bess}{\begin{eqnarray*}}
\newcommand{\eess}{\end{eqnarray*}}

\newcommand{\R}{\mathbb{R}}
\newcommand{\bpf}{{\bf Proof:\ \ }}
\newcommand{\epf}{\mbox{}\hfill $\Box$}

\begin{document}
\thispagestyle{empty}

\title{The spatial-temporal risk index and spreading dynamics for a time-periodic diffusive WNv model\thanks{The work is partially supported by the NNSF of China (Grant No. 11771381, 11701206) and CIHR and NSERC of Canada.}}
\date{\empty}
\author{Jing Ge, Zhigui Lin$\thanks{Corresponding author. Email: zglin68@hotmail.com (Z. Lin).}$, Huaiping Zhu\\
{\small  School of Mathematical Science, Yangzhou University, Yangzhou 225002, China}\\
{\small Laboratory of Mathematical Parallel Systems (LAMPS)}\\
 {\small Department of Mathematics and Statistics}\\
 {\small York University, Toronto, ON, M3J 1P3, Canada}}

 \maketitle

\begin{quote}
\noindent
{\bf Abstract.} { 
This paper is concerned with a simplified epidemic model for West Nile virus in a heterogeneous time-periodic environment. By means of the model, we will explore the impact of spatial heterogeneity of environment and temporal periodicity on the persistence and eradication of West Nile virus. The free boundary is employed to represent the
moving front of the infected region. The basic reproduction number $R_0^D$ and the spatial-temporal risk index $R_0^F(t)$, which depend on spatial heterogeneity, temporal periodicity and spatial diffusion, are defined
by considering the associated linearized eigenvalue problem. Sufficient conditions for the spreading and vanishing of West Nile virus are presented for the spatial dynamics of the virus.
 }


\end{quote}

\section{Introduction}

West Nile virus (WNv), which was first identified
in 1937 from the blood of a febrile woman in the West Nile District of Ugandan during the research on yellow fever virus \cite{ACT}, is transmitted among mosquitoes, birds, human, and other domestic animals. It is believed that WNv is long-standing in natural world in a mosquito-bird-mosquito transmission cycle \cite{CMLG}. Since the first outbreak in New York in the late summer of 1999, WNv has been spreading through the whole continent of North America for the last several years \cite{CDC2003}.  It is reported that about 1 in 5 people who are infected will develop a fever and less than $1\%$ of infected people develop a serious, sometimes fatal illness. However, there are no medications to treat or vaccines to prevent WNv infection. It is essential
to acquire some insights into the transmission dynamics of WNv in the mosquito-bird population.

There have been intensive modeling and analysis for the temporal transmission dynamics of WNv since 1999,
see for example Bowman et al. \cite{BGP}, Lewis et al. \cite{LRD}, Wan and Zhu \cite{WZ}, Abdelrazec et al. \cite{Abdelrazec2014} and so on. It is worth mentioning that Lewis et al. \cite{LRD} investigated the following simplified WNv model
\begin{eqnarray}
\left\{
\begin{array}{ll}
\frac{\partial I_b}{\partial t}=D_1\Delta I_b +\alpha_{b}\beta_{b}\frac{(N_b-I_b)}{N_b} I_m -\gamma_{b}I_b ,&(x,t)\in \Omega \times (0,+\infty), \\
\frac{\partial I_m }{\partial t}=D_2\Delta I_m +\alpha_{m}\beta_b \frac{(A_m-I_m)}{N_b} I_b -d_{m} I_m ,&(x,t)\in \Omega \times (0,+\infty),\\
I_b(x,0)=I_{b,0}(x),\ I_m(x,0)=I_{m,0}(x),&x\in \overline \Omega,
\end{array} \right.
\label{Aa1}
\end{eqnarray}
where the constants $N_b$ and $A_m$ denote, respectively, the total population of birds and adult mosquitos; $I_b(x,t)$ and $I_m(x,t)$ stand for the populations of infected birds and mosquitos at the location $x$ in the habitat $\Omega\subset \R^n$ and at time $t (\geq 0)$. The positive constants $D_1$ and $D_2$ are the diffusion coefficients for birds and mosquitoes, respectively. The
remaining parameters in the above system are described as follows:

$\bullet$ $\alpha_{m}$, $\alpha_{b}$ : WNv transmission probability per bite to mosquitoes and birds, respectively;

$\bullet$ $\beta_{b}$ : biting rate of mosquitoes on birds;

$\bullet$ $d_m $ : death rate of adult mosquitos induced by WNv;

$\bullet$ $\gamma_b $ : bird recovery rate from WNv.

In \cite{LRD}, Lewis et al. explored the spatial spread of West Nile virus, and established the existence of traveling waves
as well as computed the spatial spread speed of the infected. There are some recent studies concerning the WNv dynamics, see for example,
\cite{LZ} and references therein. However, most existing work studies the transmission of WNv in homogeneous environment and the corresponding systems are spatially-independent.

To better understand the impact of spatial diffusion and environmental heterogeneity on the transmission of infectious disease, Allen et al. \cite{AL} proposed an SIS epidemic reaction-diffusion model in a fixed domain
subject to null Neumann boundary condition
\begin{eqnarray}
\left\{
\begin{array}{lll}
S_{t}-d_S\Delta S=-\frac{\beta (x) SI}{S+I}+\gamma(x)I,\; &  x\in\Omega, \ t>0, \\
I_{t}-d_I\Delta I=\frac{\beta (x) SI}{S+I}-\gamma(x)I,\; &   x\in \Omega,\ t>0,\\
\frac {\partial S}{\partial \eta}=\frac {\partial I}{\partial \eta}=0,\;&  x\in \partial \Omega, \ t>0,
\end{array} \right.
\label{Aa11}
\end{eqnarray}
where $S(x, t)$ and $I(x, t)$ represent the susceptible and infected individuals at
location $x$ and time $t$, respectively, the positive constants $d_S$ and $d_I$ denote the corresponding diffusion rates for the susceptible and infected individuals, $\beta(x)$ and $\gamma(x)$ are positive H$\ddot{o}$lder continuous functions, which represent spatial dependent rates of disease contact transmission and disease recovery at $x$, respectively. The term $\frac{\beta (x) SI}{S+I}$ is the standard incidence of disease.
It was shown that environmental heterogeneity can influence the persistence and eradication of infectious diseases and it could cause complicated and abundant dynamics.
Recently, Peng and co-workers \cite{LPW, PZ} further investigated the asymptotical
behavior and global stability of the endemic equilibrium for system (\ref{Aa11}) subject to the Neumann boundary conditions.  In \cite{CL}, Cui and Lou considered the common effects of the diffusion and advection for an SIS epidemic model in heterogeneous environment and introduced the basic reproduction number $R_0$ for advection rate and mobility of the infected individuals.
They found that for low-risk domain, there may exist a critical value for the advection rate, under which the
disease-free equilibrium changes its stability at least twice as $d_I$ varies from zero to infinity, while the disease-free equilibrium is unstable for any $d_I$ when the advection rate is bigger than the critical value.

In most previous works, environmental heterogeneity is introduced by non-constant contact transmission and recovery rates, the related reaction-diffusion problems in a bounded domain are usually proposed to describe the persistence and eradication of infectious diseases in the fixed environment.  However, as we know, changing or expanding of an infected area is an successive process,  another remarkable feature of spatial spreading of
an infection. Mathematically, such unknown changing area is usually modeled by a free boundary problem.
Recently, there has been growing interest in understanding the free boundary and its role  in
mathematical ecology. For example, Du and Lin \cite{DL} proposed a diffusive logistic model in homogeneous environment:
\begin{eqnarray}
\left\{
\begin{array}{lll}
u_{t}-d u_{xx}=(a-bu)u,\; &0<x<h(t),\; t>0,   \\
u_x(0,t)=u(h(t),t)=0,&t>0,\\
h'(t)=-\mu u_{x}(h(t),t),  & t>0, \\
 h(0)=h_0>0,\, u(x,0)=u_{0}(x), & 0\leq x\leq h_0,
\end{array} \right.
\label{a20}
\end{eqnarray}
where the free boundary $x=h(t)$ represents the moving front of an invasive species. The spreading-vanishing dichotomy, sharp criteria for spreading and vanishing, and the asymptotic spreading speed of the free boundary problem have been established, where the asymptotic spreading speed is smaller than the minimal speed of the traveling waves of the corresponding Cauchy problem.
Since then, the study of the species invasions attracts much more attention. For the one species case, many authors explored the corresponding free boundary problems with general reaction terms $f(u)$ instead of $u(a-bu)$, such as monostable, bistable and combustion types, and obtained rather more complex description on the long time behavior of the solutions, see \cite{DLOU, DMZ, Kaneko, KM, LLZ} and reference therein. For the two species case, the competition models with free boundaries were studied in \cite{DL2, Wang16, WU}, Refs. \cite{Wang14, Wangzhao2} considered two species predator-prey models with free boundaries, and two species mutualistic model with free boundaries in a homogeneous environment was discussed in \cite{LLin}.

The spatial spreading of mosquito-borne diseases or general vector-borne diseases are much more complicated since it involves not only two species of population, but also a virus or diseases transmitted by vectors.
 Recently, it is recognized \cite{Ahn, GKLZ, HW, LZ} that the spreading of the infected environment depends on time $t$ and its fronts can be described by a free boundary. As for the impact of the spatial heterogeneity of environment in the transmission of infectious diseases, we mention the recent work \cite{GKLZ}, where they adopted a novel approach to describe the dynamical behaviors of infectious diseases. They introduced the risk index, which is related to the infected interval at time $t$, to characterize the spreading and vanishing phenomenon of infectious diseases.

Owing to the seasonal fluctuation and periodic availability of vaccination strategies and so on, the diffusion of infectious diseases varies periodically in time. The periodicity has been causing comprehensive attention in the investigation of transmission of infectious diseases. For instance, Peng and Zhao \cite{PZ} studied a reaction-diffusion SIS epidemic model in a time-periodic environment. In a recent paper \cite{GLZ1}, the authors considered a simplified SIS epidemic model with free boundaries in heterogeneous time-periodic environment.

Inspired by the former works, in present paper we will concentrate on the impact induced by spatial-temporal heterogeneity of environment in a diffusive WNv model with free boundary:
\begin{eqnarray}
\left\{
\begin{array}{ll}
\frac{\partial I_{b}}{\partial t}=D_1\frac{\partial^2I_b}{\partial x^2}-\gamma_{b}(x,t)I_{b}+\alpha_{b}\beta_{b}(x,t)\frac{(N_b -I_b )}{N_b} I_m ,&0<x<h(t),\, t>0, \\
\frac{\partial I_m}{\partial t}=D_2\frac{\partial^2I_m}{\partial x^2}-d_{m}I_{m}+\alpha_{m}\beta_{b}(x,t)\frac{(A_m -I_m)}{N_b} I_{b},&0<x<h(t),\, t>0,\\
I_b(h(t),t)=I_m(h(t),t)=0,\, &  t\geq 0,\\
\frac{\partial I_{b}}{\partial x}(0,t)=\frac{\partial I_{m}}{\partial x}(0,t)=0, &t\geq 0,\\
 h(0)=h_0, \; h'(t)=-\mu \frac{\partial I_b}{\partial x}(h(t), t), & t>0,\\
I_b(x,0)=I_{b,0}(x),\ I_m(x,0)=I_{m,0}(x),&0\leq x\leq h_0,
\end{array} \right.
\label{a3}
\end{eqnarray}
where $x=h(t)$ is the spreading front to be determined together with the infected birds $I_b(x,t)$ and infected mosquitos $I_m(x,t)$. The positive constant $\mu$ measures the expanding capability of the infected birds transmitting and diffusing towards the new area. $\beta_b(x,t), \gamma_b(x,t)\in C^{\nu_0, \frac{ \nu_0} 2}(\R \times [0,\infty))$ for some $\nu_0\in (0,1)$, which represent the biting rate of mosquitoes on birds and bird recovery rate from WNv at location $x$ and time $t$, respectively. We assume that
$\beta_b(x,t)$ and $\gamma_b(x,t)$ are positive and bounded, that is, there exist positive constants $\beta_{1}, \beta_{2}, \gamma_{1}$ and $\gamma_{2}$ such that $\beta_{1}\leq \beta_b(x,t)\leq\beta_{2}$ and $\gamma_{1}\leq\gamma_b (x,t)\leq \gamma_{2}$ in $\R \times [0,\infty)$.
Considering environmental heterogeneity, we also assume that $\beta_b(x,t), \gamma_b(x,t)$
 are periodic in $t$ with the same period $T$ (i.e., $\beta_b(x,t+T)=\beta_b(x,t)$, $\gamma_b(x,t+T)=\gamma_b(x, t)$ for all $t\geq 0$).
 The initial functions $I_{b,0}$ and $I_{m,0}$ are nonnegative and satisfy
\begin{eqnarray}
\left\{
\begin{array}{l}
I_{b,0}\in C^2([0, h_0]),\, I_{b,0}( h_0)=0,\, \frac{\partial I_{b,0}}{\partial x}(0)=0\,\textrm{and} \, 0< I_{b,0}(x) \leq N_b ,\  x\in (0, h_0), \\
I_{m,0}\in C^2([0, h_0]), I_{m,0}( h_0)=0,\, \frac{\partial I_{m,0}}{\partial x}(0)=0\,\textrm{and} \, 0< I_{m,0}(x) \leq A_m ,\  x\in (0, h_0),
\end{array} \right.
\label{Ae1}
\end{eqnarray}
where the condition (\ref{Ae1}) indicates that at initial time, the infected birds and  mosquitoes only exist in the area with $x\in (0, h_0)$,
while for the area $x \geq h_0$, no infected birds and  mosquitoes exist. Therefore,
the model means that beyond the free boundary $x=h(t)$, there is only susceptible, no infected. The equation governing the free boundary, the moving front,
$h'(t)=-\mu \frac{\partial I_b}{\partial x}(h(t),t)$ is the special situation of the well-known
Stefan condition. We notice that the similar free boundary conditions have been applied in ecological models in several earlier papers, such as in \cite{LIN, MAS2}.

The remainder of this paper is arranged as follows. In the next section, the global
existence and uniqueness of the solution to
(\ref{a3}) are presented by applying a contraction mapping theorem, and the
comparison principle is also employed. Section 3 is devoted to introducing the spatial-temporal
risk index and deriving their analytical properties, and section 4 deals with the T-periodic boundary value problem in half space. Sufficient conditions for the disease to vanish or spread and the long-time dynamical behavior are given in section 5.

\section{Preliminaries}

In this section, we first exhibit the global existence, uniqueness, regularity and some estimates on solutions of problem \eqref{a3}, we omit the proof since it is classical, and which are essentially parallel to Lemma 2.2, Theorems 2.1 and 2.2 in \cite{GLZ}.

\begin{thm} For any given $I_{b,0}(x), I_{m,0}(x)$ satisfying \eqref{Ae1}, and any $\nu \in (0, 1)$, problem \eqref{a3}
uniquely admits a global solution
$$(I_b, I_m; h)\in C^{ 1+\nu,(1+\nu)/2}(D^\infty)\times C^{ 1+\nu,(1+\nu)/2}(D^\infty)\times C^{1+\nu/2}([0,+\infty));$$
where $D^\infty:=\{(x,t)\,\vert \, x\in[0, h(t)], t\in[0,+\infty)\}$.
Moreover,
\[
0<I_b(x, t)\leq N_b, \; 0<I_m(x, t)\leq A_m\  \mbox{ for } 0<x<
h(t),\; t\in (0,+\infty),
\]
\[
 0<h'(t)\leq C, \;  \; t\in (0,+\infty), \]
for some constant $C$.\label{exist}
\end{thm}
\bpf  The local existence, uniqueness and regularity of the solution to problem \eqref{a3} can be obtained by similar methods as in Lemma 2.2, Theorems 2.1 and 2.2 in \cite{GLZ}.
We next derive the estimates of the unknown $I_b$ and $I_m$. For any given $T$, considering null Neumann boundary condition on the left boundary, we first extend the solution $(I_b, I_m; h(t))$ to $(u, v, h(t))$ such that
$$u(x,t)=I_b(|x|, t),\ v(x,t)=I_m(|x|,t)\ \textrm{for}\ -h(t)<x<h(t),\, 0\leq t\leq T,$$
then $u$ and $v$ satisfy
\begin{eqnarray}
\left\{
\begin{array}{ll}
\frac{\partial u}{\partial t}=D_1\frac{\partial^2 u}{\partial x^2}-\gamma_{b}(|x|,t)u+\alpha_{b}\beta_{b}(|x|,t)\frac{(N_b -u )}{N_b} v ,&-h(t)<x<h(t),\, 0<t\leq T, \\
\frac{\partial v}{\partial t}=D_2\frac{\partial^2 v}{\partial x^2}-d_{m}v+\alpha_{m}\beta_{b}(|x|,t)\frac{(A_m -v)}{N_b} u,&-h(t)<x<h(t),\, 0<t\leq T,\\
u(\pm h(t),t)=v(\pm h(t),t)=0,\, &  0<t\leq T,\\
u(x,0)=I_{b,0}(|x|),\ v(x,0)=I_{m,0}(|x|),&-h_0\leq x\leq h_0,
\end{array} \right.
\label{a3b}
\end{eqnarray}

 We now show that $(u, v)\geq (0, 0)$ for $-h(t)\leq x\leq h(t)$,\, $0\leq t\leq T$.
 Letting $U=u e^{-Kt}$ and $V=v e^{-Kt}$, we obtain
 \begin{eqnarray*}
\left\{
\begin{array}{ll}
\frac{\partial U}{\partial t}=D_1\frac{\partial^2 U}{\partial x^2}+(-K-\gamma_{b}-\alpha_{b}\beta_{b}\frac{v}{N_b})U+\alpha_{b}\beta_{b}V,&-h(t)<x<h(t),\, 0<t\leq T, \\
\frac{\partial V}{\partial t}=D_2\frac{\partial^2 V}{\partial x^2}+(-K-d_{m}-\alpha_{m}\beta_{b}\frac{u}{N_b} )V+\alpha_{m}\beta_{b}\frac{A_m}{N_b} U,&-h(t)<x<h(t),\, 0<t\leq T,\\
U(\pm h(t),t)=V(\pm h(t),t)=0,\, &  0<t\leq T,\\
U(x,0)=I_{b,0}(|x|),\ V(x,0)=I_{m,0}(|x|),&-h_0\leq x\leq h_0,
\end{array} \right.
\end{eqnarray*}
where $K$ is sufficiently large such that
$$K\geq 1+\alpha_{b}\beta_{b}(|x|,t)\frac{N_b+|v(x,t)|}{N_b}+\alpha_{m}\beta_{b}(|x|,t)\frac{A_m+|u(x,t)|}{N_b}$$
for $-h(t)\leq x\leq h(t)$, $0\leq t\leq T.$

We claim that $\min \{\min_{[0, T]\times [-h_0, h_0]} U, \min_{[0, T]\times [-h_0, h_0]} V\}:=m\geq 0$.
In fact, if $m<0$, then there exists $(x_0, t_0)\in\mathbb{R}^2$ with $0<t_0\le T$ and $-h(t_0)<x_0<h(t_0)$ such that $U(x_0, t_0)=m<0$, or there exists
$(x_1, t_1)\in\mathbb{R}^2$ with $0<t_1\le T$ and $-h(t_0)<x_1<h(t_0)$ such that $V(x_1, t_1)=m<0$. For the former case, $(U_{t}-D_1 U_{xx})(x_0, t_0)\leq 0$, but
\begin{align*}
[(-K-\gamma_{b}-\alpha_{b}\beta_{b}\frac{v}{N_b})U+\alpha_{b}\beta_{b}V](x_0, t_0)\geq (-K+\alpha_{b}\beta_{b}\frac{|v|}{N_b})m+\alpha_{b}\beta_{b}m\geq -m>0.
\end{align*}
For the latter case, $(V_{t}-D_2 V_{xx})(x_1, t_1)\leq 0$, but
\begin{align*}
[(-K-d_{m}-\alpha_{m}\beta_{b}\frac{u}{N_b} )V+\alpha_{m}\beta_{b}\frac{A_m}{N_b} U](x_1, t_1)\geq (-K+\alpha_{m}\beta_{b}\frac{|u|}{N_b} )m+\alpha_{m}\beta_{b}\frac{A_m}{N_b} m>0.
\end{align*}
Both are impossible. Therefore $m\geq 0$, that is $U\geq 0$, $V\geq 0$ and thus $u\geq 0, v\geq 0$
  for $ -h(t)\leq x \leq h(t)$, $0\leq t\leq T$.

 Let $(w, z)=(N_b-u, A_m-v)$, then $(w,z)$ satisfies
 \begin{eqnarray}
\left\{
\begin{array}{ll}
\frac{\partial w}{\partial t}\geq D_1\frac{\partial^2 w}{\partial x^2}-\alpha_{b}\beta_{b}(|x|,t)\frac{v}{N_b}w ,&-h(t)<x<h(t),\, 0<t\leq T, \\
\frac{\partial z}{\partial t}\geq D_2\frac{\partial^2 z}{\partial x^2}-\alpha_{m}\beta_{b}(|x|,t)\frac{u}{N_b}z,&-h(t)<x<h(t),\, 0<t\leq T,\\
w(\pm h(t),t)=N_b,\ z(\pm h(t),t)=A_m,\, &  0<t\leq T,\\
w(x,0)=N_b-I_{b,0}(|x|)\geq 0,&-h_0\leq x\leq h_0,\\
z(x,0)=A_m-I_{m,0}(|x|)\geq 0,&-h_0\leq x\leq h_0.
\end{array} \right.
\label{a3b4}
\end{eqnarray}
 Applying the maximum principle gives that $(w,z)\geq (0,0)$ for $ -h(t)\leq x \leq h(t)$, $0\leq t\leq T$.
 We then have  $(0,0)\leq (u,v)\leq (N_b, A_m)$ for $ -h(t)\leq x \leq h(t)$, $0\leq t\leq T$, which implies that $(0, 0)\leq (I_b, I_m)\leq (N_b, A_m)$ for $ 0\leq x \leq h(t)$, $0\leq t\leq T$. Moreover, using the strong maximum principle yields $(0, 0)<(I_b, I_m)$ for $ 0\leq x <h(t)$, $0<t\leq T$.

 The estimates for $h'(t)$ is followed from the maximum principle, we omit the proof since it is standard.
Noting that the bounds for $I_b, I_m$ and $h'(t)$ in $\{(x,t)\,\vert \, x\in[0, h(t)], t\in[0, T]\}$ are independent of $T$, we can use Zorn's Lemma to conclude that the solution is global and all estimates hold for $0<t<\infty$, see also Theorem 2.3 in \cite{DL}.
\epf

\bigskip
In order to facilitate later applications, we state the comparison principle, which is similar to Lemma 2.2 in \cite{GLZ1}.
\begin{lem}
$($Comparison Principle$)$
  Assume that $T\in (0,\infty)$, $\overline g, \overline h\in C^1([0,T])$, $\overline I_b (x, t)$, $\overline I_m (x,t)\in C(\overline{D}_T)\cap C^{2,1}(D_T)$ with $D_T:=\{(x,t)\vert x\in[0, h(t)], t\in (0, T]\}$ and
\begin{eqnarray*}
\left\{
\begin{array}{lll}
\frac{\partial \overline I_b}{\partial t}\geq D_1\frac{\partial^2 \overline I_b}{\partial x^2}-\gamma_{b}(x,t)\overline I_b +\alpha_{b}\beta_{b}(x,t)\frac{(N_b -\overline I_b)}{N_b} \overline I_m, &0<x<\overline h(t),\quad  0<t\leq T, \\
\frac{\partial \overline I_m}{\partial t}\geq D_2\frac{\partial^2 \overline I_m}{\partial x^2}-d_{m}\overline I_m +\alpha_{m}\beta_{b}(x,t)\frac{(A_m -\overline I_m )}{N_b} \overline I_b ,&0<x<\overline h(t),\quad  0<t\leq T,\\
\overline I_b (x,t)=\overline I_m (x, t)=0,\, &  x= \overline h(t),
\quad  0<t\leq T,\\
\frac{\partial \overline I_{b}}{\partial x}(0,t)\leq 0,\;\frac{\partial \overline I_{m}}{\partial x}(0,t)\leq 0,\;&0< t\leq T,\\
\overline h(0)\geq h_0, \; \overline h'(t)\geq -\mu \frac{\partial \overline I_b }{\partial x}(\overline h(t), t), \quad & 0<t\leq T,\\
N_b\geq \overline I_b (x,0)\geq I_{b,0}(x),\ A_m\geq\overline I_m (x,0)\geq I_{m,0}(x),&0\leq x\leq h_0.
\end{array} \right.
\end{eqnarray*}
If $\overline I_b (x,t)\leq N_b$ and $\overline I_m (x,t)\leq A_m$ in $\overline{D}_T$.
Then the solution $(I_b, I_m;  h)$ of the free boundary problem $(\ref{a3})$ satisfies
$$h(t)\leq\overline h(t),\quad 0<t\leq T,$$
$$I_b (x, t)\leq \overline I_b (x, t),\ I_m (x, t)\leq \overline I_m (x, t),\quad x\in [0, h(t)],\ 0\leq t\leq T.$$\label{Com}
\end{lem}

It is worth mentioning that the functions in problem (\ref{a3}) are quasi-monotone nondecreasing and the system is cooperative if $\overline I_b (x,t)\leq N_b$ and $\overline I_m (x,t)\leq A_m$ in $\overline{D}_T$. Certainly we also need the conditions $I_b (x,t)\leq N_b$ and $I_m (x,t)\leq A_m$ in $\overline {D}_T$, which has been given in Theorem \ref{exist}. Biologically, it is natural since that $N_b$ is the total number of birds and $A_m$ is the total number of mosquitoes.

The pair $(\overline I_b, \overline I_m; \overline h)$ in Lemma \ref{Com} is usually called an upper solution
of problem \eqref{a3}. Similarly, we can define the lower solution $(\underline I_b, \underline I_m;  \underline h)$ of problem \eqref{a3} by reversing all the inequalities in the obvious places.

\section{The spatial-temporal risk index}

 The basic reproduction number $R_0$ is one of the most important concepts
in epidemiology, it has commonly been used to evaluate the probability of epidemics and to measure the effort needed to control an infectious disease.
 $R_0$ is defined as the expected number of secondary cases produced, in a completely susceptible population, by a typical infected individual during its entire period of infectiousness \cite{DHM1990}. For spatially-independent epidemic models, which are described by ordinary differential systems, the numbers are usually calculated by the next generation matrix method \cite{vanden2002}, while for the models constructed by reaction-diffusion systems, the numbers are formulated as the spectral radius of next infection operator induced by a new infection rate matrix and an evolution operator of an infective distribution \cite{WZ12}, and the numbers could be expressed in the term of the principal eigenvalues of relevant eigenvalue problems \cite{AL, Zhao2017}.

In this section, we first present the basic reproduction number and its properties
for the corresponding system in $[0,b)$ with $b>0$.
The basic reproduction numbers are related to the following linear periodic-parabolic eigenvalue problem:
 \be \left\{
\begin{array}{ll}
\frac{\partial \phi}{\partial t}-D_1\Delta \phi=\alpha_{b}\beta_{b}(x,t)\frac 1 R\psi-\gamma_{b}(x,t)\phi +\mu \phi,&(0, b)\times (0, +\infty), \\
\frac{\partial \psi}{\partial t}-D_2\Delta \psi=\alpha_{m}\beta_{b}(x,t)\frac{A_m }{N_b R} \phi-d_{m}\psi+\mu \psi,&(0, b)\times (0, +\infty),\\
\phi_x(0,t)=\psi_x(0,t)=\phi(b,t)=\psi(b,t)=0,\, & [0, +\infty),\\
\phi(x,t+T)=\phi(x,t),\, \psi(x,t+T)=\psi(x,t),&[0,b]\times [0, +\infty),
\end{array} \right.
\label{a31}
\ee
where $R>0$. Setting
 \be
 \mathcal{L}_R:=
 \left(
\begin{array}{ll}
\partial_t-D_1\Delta +\gamma_{b}(x,t)& -\alpha_{b}\beta_{b}(x,t)\frac 1R \\
-\alpha_{m}\beta_{b}(x,t)\frac{A_m }{N_b R} & \partial_t-D_2\Delta +d_m
\end{array} \right),
\label{a32}
\ee
then problem \eqref{a31} can be formulated as an abstract eigenvalue problem
\be
 \mathcal{L}_R
 {\phi \choose \psi}=\mu {\phi \choose \psi},
\label{a33}
\ee
in the space
$$X:=\{(\phi,\psi)\in (C^{\nu,\nu/2}([0,b]\times [0, +\infty)))^2\,:\phi, \psi\, \textrm{are}\, T-\textrm{periodic in}\, t\},$$
and the domain of the operator dom$(\mathcal{L}_R)=X_1$ is defined by
$$X_1=\{(\phi,\psi)\in (C^{2+\nu, 1+\nu/2}([0,b]\times [0, +\infty)))^2\,: \phi_x(0,t)=\psi_x(0,t)=\phi(b,t)=\psi(b,t)=0$$
$$\qquad \qquad \textrm{for}\, t\in [0, +\infty),\ \phi, \psi\, \textrm{are}\, T-\textrm{periodic in}\, t\}.$$

For any given $R>0$, system \eqref{a33} is strongly cooperative in the sense that $\alpha_{b}\beta_{b}(x,t)>0$ and $\alpha_{m}\beta_b(x,t) \frac{A_m}{N_b }>0$ for all $(x,t)\in [0,b]\times [0, +\infty)$.
Similarly as in \cite{ADP, ALG}, it follows from the Krein-Rutman theorem ( see, e.g., Theorem 7.2 in \cite{Hess}) that there exists a unique value $\mu:=\mu_1(R, [0,b))$, and called the principal eigenvalue, such that problem \eqref{a31}, and equivalently \eqref{a33}, admits a unique solution pair $(\phi_R, \psi_R)$ (subject to constant multiples) with
$\phi_R>0$ and $\psi_R>0$ in $[0, b)\times [0, +\infty)$. The solution pair $(\phi_R, \psi_R) \in X_1$ is called the principal eigenfunction corresponding to $\mu_1$.
Moreover, one can deduce from \cite{ADP, ALG} the following continuity and monotonicity.
\begin{lem} $\mu_1(R, [0,b))$ is continuous and strictly increasing with respect to $R$, and $\mu_1(R, [0, b))$ is decreasing with respect to $b$ in the sense that
$\mu_1(R, [0,b_1))>\mu_1(R, [0, b_2))$ if $b_1<b_2$.
\label{B1f}
\end{lem}

Let $R_0^D:=R^D_0([0, b))$ be the unique principal eigenvalue of the periodic-parabolic eigenvalue problem with $\mu=0$ for problem \eqref{a31},
\be \left\{
\begin{array}{ll}
\frac{\partial \phi}{\partial t}-D_1\Delta \phi=\alpha_{b}\beta_{b}(x,t)\frac 1 {R^D_0}\psi-\gamma_{b}(x,t)\phi,&(0, b)\times (0, +\infty), \\
\frac{\partial \psi}{\partial t}-D_2\Delta \psi=\alpha_{m}\beta_{b}(x,t)\frac{A_m }{N_b {R^D_0}} \phi-d_{m}\psi,&(0, b)\times (0, +\infty),\\
\phi_x(0,t)=\psi_x(0,t)=\phi(b,t)=\psi(b,t)=0,\, & [0, +\infty),\\
\phi(x,t+T)=\phi(x,t),\, \psi(x,t+T)=\psi(x,t),&[0,b]\times [0, +\infty).
\end{array} \right.
\label{a311}
\ee
The principal eigenvalue $R_0^D$ is the only positive eigenvalue admitting a unique positive eigenfunction $(\phi,\psi)$ (subject to a constant multiple).
It was proved in \cite{Zhao2017} that
$R_0^D$ is the spectral radius of the next generation operator induced by a new infection rate matrix and an evolution operator of an infective distribution.
With the above definition, we have the following relation between the two eigenvalues, see also Lemma 3.1 in \cite{GKLZ} and Theorem 11.3 in \cite{Zhao2017}.
\begin{thm}\label{r0} $1-R_0^D$ has the same sign as $\lambda_0$, where $\lambda_0:=\lambda_0([0,b))$ is the principal eigenvalue of the eigenvalue problem
\be \left\{
\begin{array}{ll}
\frac{\partial \phi}{\partial t}-D_1\Delta \phi=\alpha_{b}\beta_{b}(x,t)\psi-\gamma_{b}(x,t)\phi +\lambda_0 \phi,&(0, b)\times (0, +\infty), \\
\frac{\partial \psi}{\partial t}-D_2\Delta \psi=\alpha_{m}\beta_{b}(x,t)\frac{A_m }{N_b } \phi-d_{m}\psi+\lambda_0 \psi,&(0, b)\times (0, +\infty),\\
\phi_x(0,t)=\psi_x(0,t)=\phi(b,t)=\psi(b,t)=0,\, & [0, +\infty),\\
\phi(x,t+T)=\phi(x,t),\, \psi(x,t+T)=\psi(x,t),&[0,b]\times [0, +\infty).
\end{array} \right.
\label{a312}
\ee
\end{thm}
\bpf Comparing \eqref{a31} with \eqref{a312}, we can derive that $\lambda_0([0, b))=\mu_1(1, [0,b))$.
On the other hand, one can easily deduce from the monotonicity with respect to the coefficients in \eqref{a31} that $\lim_{R\to 0^+}\mu_1(R, [0, b))<0$ and $\lim_{R\to +\infty}\mu_1(R, [0, b))>0$,
therefore $R^D_0([0, b))$ is the unique positive root of the equation $\mu_1(R, [0, b))=0$.
Owing to $\lambda_0=\mu_1(1, [0, b))-\mu_1(R^D_0, [0, b))$, the result $sign \{1-R_0^D\}=sign \{\lambda_0\}$ follows directly from the monotonicity of $\mu_1(R, [0, b))$ with respect to $R$.
\epf
\smallskip

If all coefficients in problem \eqref{a311} are constant, we can provide an explicit formula for $R_0^D([0, b))$, which is known as the basic reproduction number for the corresponding diffusive WNv model.
\begin{thm} If $\beta_b(x,t)=\beta^*_b$, $\gamma_b(x,t)=\gamma^*_b$, then the principal eigenvalue $R^D_0$ for \eqref{a311}, or the basic reproduction number for model \eqref{a3}, is represented by
\begin{equation}\label{bac2}
R_0^D([0, b))=\sqrt{\frac{ A_m\alpha_b\alpha_m(\beta^*_{b})^2}{ N_b[D_1(\frac \pi{2b})^2+\gamma^*_b][D_2 (\frac\pi{2b})^2+d_m]}}\, .
\end{equation}
\label{basic2}
\end{thm}
\bpf
Let $$\psi^*(x)=\cos (\frac \pi{2b}x), \quad x\in [0, b],$$
$$R^*=\frac{ A_m\alpha_b\alpha_m(\beta^*_{b})^2}{ N_b[D_1(\frac \pi{2b})^2+\gamma^*_b][D_2 (\frac\pi{2b})^2+d_m]},$$
$$\phi^*(x)=\frac{ \alpha_b \beta^*_b}{\sqrt{R^*}[D_1(\frac \pi{2b})^2+\gamma^*_b]}\psi^*(x).$$
Then we know that $(\phi^*,\psi^*)$ is a positive solution of problem (\ref{a311}) with $R_0^D=\sqrt{R^*}$,
and \eqref{bac2} follows directly from the uniqueness of the principal eigenvalue of (\ref{a311}).
\epf

\bigskip

It is well-known that the basic reproduction number is a critical threshold to determine whether the disease is persistent or extinct. When we consider the spreading or vanishing phenomenon of the disease, it is often the constant defined for a spatially-independent model or a diffusive epidemic model in a fixed
region. However, for our model \eqref{a3}, the infected interval is changing with time $t$, therefore, the basic reproduction number is not a constant and
should be a function of $t$. So we here call it {\bf the spatial-temporal risk index}, which is expressed by
\begin{equation}\label{free}
R_0^F(t):=R_0^D([0, h(t))),
\end{equation}
where $R_0^D$ is the principal eigenvalue of the corresponding problem (\ref{a311}) in $[0, h(t))$.
With the above definition, we have the following properties of $R_0^F (t)$.
\begin{lem} The following statements are valid:\label{proper}
\item[$(i)$] $R_0^F(t)$ is strictly monotone increasing function with respect to $t$, that is, if $0\leq t_1<t_2$, then $R_0^F(t_1)<R_0^F(t_2)$;

\item[$(ii)$] if $h(t)\to \infty$ as $t\to \infty$, then
$$R_0^F(t)\to R_0:=\sqrt{\frac{ A_m\alpha_b\alpha_m(\beta^*_b)^2}{ N_b\gamma^*_bd_m}}\ \textrm{as}\ t\to \infty$$
provided that $\beta_b(x,t)=\beta^*_b$, $\gamma_b(x,t)=\gamma^*_b$, where $R_0$ is the usual basic reproduction number for the corresponding spatially-independent model.
\end{lem}

\section{The T-periodic boundary value problem in half line}

In order to discuss the long-time dynamical behavior of solution when spreading occurs, in what follows,
we will explore a stationary problem: the T-periodic boundary value problem in half space. The
T-periodic boundary value problem associated with the free boundary problem \eqref{a3} in half line is
\begin{eqnarray}
\left\{
\begin{array}{ll}
\frac{\partial U}{\partial t}=D_1\frac{\partial^2 U}{\partial x^2}-\gamma_{b}(x,t)U +\alpha_{b}\beta_{b}(x,t)\frac{(N_b -U )}{N_b} V ,&x>0,\,0\leq t\leq T, \\
\frac{\partial V}{\partial t}=D_2\frac{\partial^2V}{\partial x^2}-d_{m}V+\alpha_{m}\beta_{b}(x,t)\frac{(A_m -V)}{N_b} U,&x>0,\, 0\leq t\leq T,\\
\frac{\partial U}{\partial x}(0,t)=\frac{\partial V}{\partial x}(0,t)=0, &0\leq t\leq T,\\
U(x,0)=U(x,T),\ V(x,0)=V(x,T),&x\geq 0
\end{array} \right.
\label{b1}
\end{eqnarray}
which is related to the T-periodic boundary value problems in a bounded interval $(0,l)$
\begin{eqnarray}
\left\{
\begin{array}{ll}
\frac{\partial U}{\partial t}=D_1\frac{\partial^2 U}{\partial x^2}-\gamma_{b}(x,t)U +\alpha_{b}\beta_{b}(x,t)\frac{(N_b -U )}{N_b} V ,&0<x<l,\, 0\leq t\leq T, \\
\frac{\partial V}{\partial t}=D_2\frac{\partial^2V}{\partial x^2}-d_{m}V+\alpha_{m}\beta_{b}(x,t)\frac{(A_m -V)}{N_b} U,&0<x<l,\ 0\leq t\leq T,\\
(U_x, V_x)(0,t)=(0,0),\ (U, V)(l,t)=(0,0),&0\leq t\leq T,\\
U(x,0)=U(x,T),\ V(x,0)=V(x,T),&0\leq x\leq l,
\end{array} \right.
\label{b2}
\end{eqnarray}
and
\begin{eqnarray}
\left\{
\begin{array}{ll}
\frac{\partial U}{\partial t}=D_1\frac{\partial^2 U}{\partial x^2}-\gamma_{b}(x,t)U +\alpha_{b}\beta_{b}(x,t)\frac{(N_b -U )}{N_b} V ,&0<x<l,\, 0\leq t\leq T, \\
\frac{\partial V}{\partial t}=D_2\frac{\partial^2V}{\partial x^2}-d_{m}V+\alpha_{m}\beta_{b}(x,t)\frac{(A_m -V)}{N_b} U,&0<x<l,\ 0\leq t\leq T,\\
(U_x, V_x)(0,t)=(0,0),\
(U, V)(l,t)=(N_b, A_m),&0\leq t\leq T,\\
U(x,0)=U(x,T),\ V(x,0)=V(x,T),&0\leq x\leq l.
\end{array} \right.
\label{b21}
\end{eqnarray}
The boundary conditions on $x=l$ for \eqref{b2} and \eqref{b21} are different. Problem \eqref{b2} is used to construct the minimal solution of problem \eqref{b1} and problem \eqref{b21} is used to construct the maximal solution of problem \eqref{b1}.
To study problems \eqref{b1}, \eqref{b2} and \eqref{b21}, we need to consider the corresponding initial boundary problem to \eqref{b1} in half space
\begin{eqnarray}
\left\{
\begin{array}{ll}
\frac{\partial u}{\partial t}=D_1\frac{\partial^2 u}{\partial x^2}-\gamma_{b}(x,t)u +\alpha_{b}\beta_{b}(x,t)\frac{(N_b -u )}{N_b} v ,&x>0,\, t>0, \\
\frac{\partial v}{\partial t}=D_2\frac{\partial^2v}{\partial x^2}-d_{m}v+\alpha_{m}\beta_{b}(x,t)\frac{(A_m -v)}{N_b} u,&x>0,\, t>0,\\
\frac{\partial u}{\partial x}(0,t)=\frac{\partial v}{\partial x}(0,t)=0, &t>0,\\
u(x,0)=I_{b,0}(x), \ v(x,0)=I_{m,0}(x),&x\geq 0
\end{array} \right.
\label{b3}
\end{eqnarray}
where $I_{b,0}(x), I_{m,0}(x)$ are non-trivial continuous functions and satisfy $(0, 0)\leq (I_{b,0}, I_{m,0})\leq (N_b(x), A_m(x))$ for $x\geq 0$. We first give the estimates for solutions to problems \eqref{b2} and \eqref{b3}, which can be derived by the comparison principle.

\begin{lem} Any bounded nonnegative nontrivial solution $(U, V)$ of T-periodic boundary value problem \eqref{b1} satisfies
$$(0,0)<(U(x,t), V(x,t))<(N_b, A_m),\quad  x>0,\,0\leq t\leq T,$$
and the unique bounded solution $(u, v)$ of initial boundary problem \eqref{b3} satisfies
$$(0,0)\leq (u(x,t), v(x,t))\leq (N_b, A_m),\quad  x>0,\,t\geq 0.$$
\end{lem}

Next results present the relations of the solutions to the above problems.
\begin{lem} For any $l> L_0$, where $L_0$ satisfy $R_0^D([0,L_0))=1,$ the T-periodic boundary value problem \eqref{b2} admits the minimal positive solution $(\underline{U}_l,\underline{V}_l)$. Moreover,
the solution $(u, v)$ of problem \eqref{b3} satisfies
\begin{eqnarray}
\left\{
\begin{array}{ll}
\underline{U}_l(x,t)\leq\liminf\limits_{n\longrightarrow\infty}u(x, t+nT)\leq N_b,\\
\underline{V}_l(x,t)\leq\liminf\limits_{n\longrightarrow\infty}v(x, t+nT)\leq A_m,
\end{array} \right.
\label{ba4}
\end{eqnarray}
on $[0, l]\times[0, \infty)$.\label{ell1}
\end{lem}
\bpf Owing to $R_0^D([0,l))>1$ for any $l>L_0$, therefore, the periodic-parabolic problem
\begin{eqnarray}
\left\{
\begin{array}{ll}
\dfrac{\partial \phi}{\partial t}=D_1\dfrac{\partial^2 \phi}{\partial x^2}-\gamma_{b}(x,t)\phi +\alpha_{b}\beta_{b}(x,t) \psi +\lambda \phi,&0<x<l,\, 0\leq t\leq T, \\
\dfrac{\partial \psi}{\partial t}=D_2\dfrac{\partial^2 \psi}{\partial x^2}-d_{m}\psi+\alpha_{m}\beta_{b}(x,t)\dfrac{A_m }{N_b} \phi+\lambda\psi,&0<x<l,\ 0\leq t\leq T,\\
\phi_x(0,t)=\psi_x(0,t)=\phi (l,t)=\psi (l,t)=0, &0\leq t\leq T,\\
\phi(x,0)=\phi(x,T),\ \psi(x,0)=\psi(x,T),&0\leq x\leq l.
\end{array} \right.
\label{ba5}
\end{eqnarray}
admits the principal eigenvalue $\lambda_0(<0)$ and the corresponding eigenfunction $(\phi(x,t), \psi(x,t))$
satisfying $(\phi(x,t), \psi(x,t))>(0,0)$ in $[0,l)\times[0,T]$. It is easy to verify that, for sufficiently small $\delta$, $(N_b, A_m)$ and $(\delta\phi, \delta\psi)$ are a pair of ordered upper and lower solutions of \eqref{b2}.

Let
$$K_1=\sup_{[0, \infty)\times [0,T]} \gamma_b(x,t)+\alpha_b \frac {A_m}{N_b}\sup_{[0, \infty)\times [0,T]} \gamma_b(x,t),\
K_2=d_m+\alpha_b \sup_{[0, \infty)\times [0,T]} \gamma_b(x,t),$$
then the equations in \eqref{b2} become
\begin{eqnarray*}
\left\{
\begin{array}{ll}
\frac{\partial U}{\partial t}-D_1\frac{\partial^2 U}{\partial x^2}+K_1 U=K_1U-\gamma_{b}(x,t)U +\alpha_{b}\beta_{b}(x,t)\frac{(N_b -U )}{N_b} V:=f_1(U, V),& \\
\frac{\partial V}{\partial t}-D_2\frac{\partial^2V}{\partial x^2}+K_2 V=K_2V-d_{m}V+\alpha_{m}\beta_{b}(x,t)\frac{(A_m -V)}{N_b} U:=f_2(U,V),&
\end{array} \right.
\end{eqnarray*}
It is easy to see that $f_1$ and $f_2$ are increasing with respect to $U$ and $V$ if $(0,0)\leq (U, V)\leq (N_b, A_m)$.

Using $(\underline{U}^{(0)},\underline{V}^{(0)})=(\delta\phi, \delta\psi)$
as initial iteration, we construct a sequence ${(\underline{U}^{(n)},\underline{V}^{(n)})}$ from
the linear boundary problem
\begin{eqnarray}
\left\{
\begin{array}{ll}
\frac{\partial U^{(n)}}{\partial t}-D_1\frac{\partial^2 U^{(n)}}{\partial x^2}+K_1 U^{(n)}=f_1(U^{(n-1)}, V^{(n-1)}), &0<x<l,\, t>0, \\
\frac{\partial V^{(n)}}{\partial t}-D_2\frac{\partial^2V^{(n)}}{\partial x^2}+K_2 V^{(n)}=f_2(U^{(n-1)}, V^{(n-1)}),&0<x<l,\, t>0,\\
U^{(n)}_x(0,t)=V^{(n)}_x(0,t)=U^{(n)}(l,t)=V^{(n)}(l,t)=0, &t>0,\\
U^{(n)}(x,0)=U^{(n-1)}(x,T),\ V^{(n)}(x,0)=V^{(n-1)}(x,T),&0\leq x\leq l.
\end{array} \right.
\label{b22}
\end{eqnarray}
Moreover, it follows from monotonicity of $f_1$ and $f_2$ that the well-defined sequences  ${(\underline{U}^{(n)},\underline{V}^{(n)})}$ possess the monotone property
$$(\delta\phi, \delta\psi)\leq (\underline{U}^{(n)},\underline{V}^{(n)})\leq (\underline{U}^{(n+1)},\underline{V}^{(n+1)})\leq (N_b, A_m)$$
in $[0, l]\times [0, +\infty)$ for every $n=1,2, \cdots.$
Therefore, the limits of the sequences
 $$\lim\limits_{n\rightarrow \infty}{(\underline{U}^{(n)},\underline{V}^{(n)})}=(\underline{U}_l,\underline{V}_l)$$ exist and the limit $(\underline{U}_l,\underline{V}_l)$ is a solution of \eqref{b2}.

 We now claim that it is also the minimal positive solution of \eqref{b2}. If fact, for any positive solution $({U}_l, {V}_l)$, for small $\delta$, $({U}_l, {V}_l)$ and $(\delta\phi, \delta\psi)$ are a pair of ordered upper and lower solutions of \eqref{b2}. By the same iteration given by \eqref{b22}, we can derive that
 $$(\delta\phi, \delta\psi)\leq (\underline{U}^{(n)},\underline{V}^{(n)})\leq (\underline{U}^{(n+1)},\underline{V}^{(n+1)})\leq ({U}_l, {V}_l)$$
and then $(\underline{U}_l,\underline{V}_l)\leq ({U}_l, {V}_l)$ in $[0, l]\times [0, +\infty)$.

Next, let $(u, v)$ be the solution of problem \eqref{b3} for $(x,t)\in [0, \infty)\times [0, \infty)$ with nontrivial nonnegative initial value, then
$$(0, 0)<(u, v)(x,T)<(N_b, A_m), \quad x\in [0, \infty)$$
and there exists $\delta>0$ such that
$$(\delta \phi, \delta \psi)(x,0)\leq (u, v)(x, T)\leq(N_b, A_m), \quad x\in [0, l].$$
Consider the system \eqref{b3} with the initial condition $(u, v)(x,T)$ in $[0, l]$. Since by the initial condition in \eqref{b22} for $n=1$, $U^{(1)}(x,0)=U^{(0)}(x,T)=\delta \phi(x,T)=\delta \phi(x,0)$, $V^{(1)}(x, 0)=V^{(0)}(x,T)=\delta \psi(x, T)=\delta \psi(x,0)$ in $[0, l]$. By comparison principle, we see that
$$(\underline{U}^{(1)},\underline{V}^{(1)} )(x,t)\leq (u, v)(x,t+T)\leq (N_b, A_m)$$ on $[0, l]\times [0, \infty)$.
Similarly as Lemma 3.2 in \cite{P6}, by using the comparison principle and the principle of induction, we have that
$$(\underline{U}^{(n)},\underline{V}^{(n)} )(x,t)\leq (u, v)(x,t+nT)\leq (N_b, A_m)$$ on $[0, l]\times [0, \infty)$ for every $n=1,2,\cdots$, which concludes the desired result \eqref{ba4}.
 \epf
 \smallskip

\begin{lem} For any $l> L_0$, where $L_0$ satisfy $R_0^D([0,L_0))=1,$ the T-periodic boundary value problem \eqref{b21} admits the maximal positive solution $(\overline{U}_l,\overline{V}_l)$. Moreover,
the solution $(u, v)$ of problem \eqref{b3} satisfies
\begin{eqnarray}
\left\{
\begin{array}{ll}
\underline{U}_l(x,t)\leq\limsup\limits_{n\longrightarrow\infty}u(x,t+nT)\leq \overline{U}_l(x,t),\\
\underline{V}_l(x,t)\leq\limsup\limits_{n\longrightarrow\infty}v(x, t+nT)\leq \overline{V}_l(x,t)
\end{array} \right.
\label{ba4-1}
\end{eqnarray}
on $[0, l]\times[0, \infty)$.\label{le43}
\end{lem}
\bpf The proof is similar to that of Lemma \ref{ell1}, we give the sketch here. First, we can see that $(N_b, A_m)$ and $(\underline{U}_l, \underline{V}_l)$ are a pair of ordered upper and lower solutions of \eqref{b21}.

Using $(\overline{U}^{(0)},\overline{V}^{(0)})=(N_b, A_m)$
as initial iteration, we construct a sequence ${(\overline{U}^{(n)},\overline{V}^{(n)})}$ from
the linear boundary problem
\begin{eqnarray}
\left\{
\begin{array}{ll}
\frac{\partial U^{(n)}}{\partial t}-D_1\frac{\partial^2 U^{(n)}}{\partial x^2}+K_1 U^{(n)}=f_1(U^{(n-1)}, V^{(n-1)}), &0<x<l,\, t>0, \\
\frac{\partial V^{(n)}}{\partial t}-D_2\frac{\partial^2V^{(n)}}{\partial x^2}+K_2 V^{(n)}=f_2(U^{(n-1)}, V^{(n-1)}),&0<x<l,\, t>0,\\
(U^{(n)}_x, V^{(n)}_x)(0,t)=0,\ (U^{(n)},V^{(n)})(l,t)=(N_b, A_m), &t>0,\\
U^{(n)}(x,0)=U^{(n-1)}(x,T),\ V^{(n)}(x,0)=V^{(n-1)}(x,T),&0\leq x\leq l.
\end{array} \right.
\label{b2-2}
\end{eqnarray}
Moreover, the well-defined sequences  ${(\overline{U}^{(n)},\overline{V}^{(n)})}$ possess the monotone property
$$(\underline{U}_l, \underline{V}_l)\leq (\overline{U}^{(n+1)},\overline{V}^{(n+1)})\leq (\overline{U}^{(n)},\overline{V}^{(n)})\leq (N_b, A_m)$$
in $[0, l]\times [0, +\infty)$ for every $n=1,2, \cdots.$
Therefore, the limits of the sequences
 $$\lim\limits_{n\rightarrow \infty}{(\overline{U}^{(n)},\overline{V}^{(n)})}=(\overline{U}_l,\overline{V}_l)$$ exist and the limit $(\overline{U}_l,\overline{V}_l)$ is a solution of \eqref{b21}.

 We now claim that it is also the maximal solution of \eqref{b2}. In fact, for any positive solution $({U}_l, {V}_l)$, $(N_b, A_m)$ and $({U}_l, {V}_l)$
 are a pair of ordered upper and lower solutions of \eqref{b2}. By the same iterative procedure given in \eqref{b2-2}, we can derive that
 $$(U_l, V_l)\leq (\overline{U}^{(n+1)},\overline{V}^{(n+1)})\leq (\overline{U}^{(n)}, \overline{V}^{(n)})\leq (N_b, A_m)$$
and then $({U}_l, {V}_l)\leq (\overline{U}_l, \overline {V}_l)$ in $[0, l]\times [0, \infty)$.

Similarly as Lemma 3.2 in \cite{P6}, by using the comparison principle and the principle of induction, we have that
$$(u, v)(x, t+nT)\leq (\overline{U}^{(n)}, \overline{V}^{(n)} )(x, t)\leq (N_b, A_m)$$ on $[0, l]\times [0, \infty)$ for every $n=1,2,\cdots$, which concludes the desired result \eqref{ba4-1}.
 \epf

\begin{thm}\label{theb5}
Suppose that $R^D_0([0, +\infty))>1$ hold. Then T-periodic boundary value problem \eqref{b1} admits the maximal and the minimal positive periodic solutions $(\overline{U}, \overline{V})$ and
$(\underline{U},\underline{V})$. Moreover,
\begin{eqnarray}
\left\{
\begin{array}{ll}
\underline{U}(x,t)\leq\liminf\limits_{n\longrightarrow\infty}u(x, t+nT)\leq
\limsup\limits_{n\longrightarrow\infty}u(x,t+nT)\leq\overline{U}(x,t),\\
\underline{V}(x,t)\leq\liminf\limits_{n\longrightarrow\infty}v(x,t+nT)\leq
\limsup\limits_{n\longrightarrow\infty}v(x, t+nT)\leq\overline{V}(x,t)
\end{array} \right.
\label{b5}
\end{eqnarray}
locally uniformly in $[0, \infty)\times[0,T]$, where $(u,v)$ is the unique solution of problem \eqref{b3}.
\end{thm}
\bpf
The proof is based on the upper and lower solutions methods and will be divided into three steps.

Step 1. The construction of $(\overline{U},\overline{V})$ and $(\underline{U}, \underline{V})$.

Owing to the assumption that $R_0^D([0, +\infty))>1$, there exists a unique $L_0$ such that $R_0^D([0, L_0))=1$.
We first present the monotonicity and show that if $L_0< l_1<l_2$, then $(\overline{U}_{l_1}, \overline{V}_{l_1})\geq (\overline{U}_{l_2}, \overline{V}_{l_2})$
and $(\underline{U}_{l_1}, \underline{V}_{l_1})\leq (\underline{U}_{l_2}, \underline{V}_{l_2})$ in $[0, l_1]$. The result is derived by comparing the boundary condition and initial conditions in \eqref{b2} and \eqref{b21} for $l=l_1$ and $l=l_2$.

Since $\overline{U}_l$ and $\overline{V}_l$ is monotone decreasing with respect to $l$, we can use the regularity theory for parabolic equations and compactness argument to deduce that $(\overline{U}_l, \overline{V}_l)(x,t)$ converge to $(\overline{U},\overline{V})(x,t)$ as $l\rightarrow+\infty$ for $(x,t)\in[0,+\infty)\times [0,T]$ and
$(\overline{U},\overline{V})(x,t)$ is a solution to the T-periodic boundary value problem \eqref{b1}. $(\underline{U}, \underline{V})(x,t)$ can be constructed by the similar way.

Step 2. We claim that $(\overline{U},\overline{V})(x,t)$ is the maximal solution to the T-periodic boundary value problem \eqref{b1}.

In fact, for any positive solution $(U, V)$ of problem \eqref{b1}, $(N_b, A_m)$ and $(U, V)$
 are a pair of ordered upper and lower solutions of \eqref{b21} in $[0, l]\times [0, \infty)$ for any $l>L_0$. By the iterative procedure given in \eqref{b2-2}, we can derive that
 $$(U, V)\leq (\overline{U}^{(n+1)},\overline{V}^{(n+1)})\leq (\overline{U}^{(n)}, \overline{V}^{(n)})\leq (N_b, A_m)$$
and then $(U, V)\leq (\overline{U}_l, \overline {V}_l)$ in $[0, l]\times [0, \infty)$, which gives that
$(U, V)\leq (\overline{U}, \overline {V})$ in $[0, +\infty)\times [0, \infty)$.
Similarly, we can prove that $(\underline{U},\underline{V})(x,t)$ is the minimal solution to the T-periodic boundary value problem \eqref{b1}.

Step 3. The proof of \eqref{b5}.

Recalling that
$$\lim_{l\rightarrow \infty}(\overline{U}_l, \overline V_l)(x,t)=(\overline{U}, \overline V)(x,t)$$
locally uniformly for $(x,t)\in[0, +\infty)\times[0,T]$, we then have, for any given $M\in (0, +\infty)$, $(\overline{U}_l, \overline {V}_l)$ converges to $(\overline{U}, \overline{V})$ uniformly on $[0, M]\times[0,T]$.
Hence, for any $\varepsilon>0$, there exists a positive constant $l_\varepsilon>M$
such that
\begin{equation}\label{ineuq1}
(\overline{U}_{l_\varepsilon}, \overline {V}_{l_\varepsilon})\leq (\overline{U}+\varepsilon, \overline{V}+\varepsilon)
\end{equation}
 on $[0, M]\times [0,T]$.

On the other hand, for the above $l_\varepsilon$, it follows from Lemma \ref{le43} that
$$\limsup\limits_{n\longrightarrow\infty} (u, v)(x,t+nT)\leq (\overline{U}_{l_\varepsilon}, \overline V_{l_\varepsilon})(x,t)$$
for $(x,t)\in [0, l_\varepsilon]\times [0,T]$. Therefore, we deduce that
$$\limsup\limits_{n\longrightarrow\infty} (u, v)(x,t+nT)\leq (\overline{U}+\varepsilon, \overline V +\varepsilon)(x,t)$$
for $(x,t)\in [0, M]\times [0,T]$. In view of the arbitrariness of $\varepsilon$, we conclude that
$$\limsup\limits_{n\longrightarrow\infty} (u, v)(x,t+nT)\leq (\overline{U}, \overline V)(x,t)$$ uniformly for $(x,t)\in [0, M]\times [0,T]$. Similarly, we have
$$\liminf\limits_{n\longrightarrow\infty} (u, v)(x,t+nT)\geq (\underline{U}, \underline V)(x,t)$$ uniformly for $(x,t)\in [0, M]\times [0,T]$.
\epf

\section{Spreading and vanishing}
In this section, some sufficient conditions for spreading or vanishing are established, as well as the long-time
dynamical behavior is presented when the spreading scenario happens.

It follows from Theorem \ref{exist} that the infected region is expanding as time $t$ increasing. In the sense that the moving front $x=h(t)$ is monotonic increasing, so
there exist $h_\infty\in (h_0, +\infty]$ such that $\lim\limits_{t\to +\infty}\ h(t)=h_\infty$.
Epidemically, it is well-known that if the infected region is bounded and the infected individuals will
die out gradually, we say the disease is vanishing, which means that the epidemic can be controlled. Mathematically, we first exhibit the following definitions.
\begin{defi}
The virus is {\bf vanishing} if
$$h_\infty <\infty\ \textrm{ and}\
 \lim_{t\to +\infty} \ (||I_b (\cdot, t)||_{C([0,h(t)])}+||I_m (\cdot, t)||_{C([0, h(t)])})=0,$$
  and  {\bf spreading} if $$h_\infty =\infty\ \textrm{and}\
\limsup_{t\to +\infty}\ (||I_b (\cdot, t)||_{C([0,h(t)])}+||I_m (\cdot, t)||_{C([0,h(t)])})>0.$$
\end{defi}

In what follows, we will theoretically present the sufficient conditions for the vanishing scenario of WNv, which can provide some effective measures and strategies for the public health administration to control West Nile virus timely.
\begin{lem}  If $h_\infty<\infty$, then we have
$$\lim_{t\to
+\infty} \ ||I_b(\cdot, t)||_{C([0, h(t)])}=\lim_{t\to
+\infty} \ ||I_m(\cdot, t)||_{C([0, h(t)])}=0.$$
\label{vanishlm}
\end{lem}
\bpf Arguing indirectly, we assume that $\displaystyle \limsup_{t\to
+\infty} \ ||I_b(\cdot, t)||_{C([0, h(t)])}=\delta>0$ for contradiction. Therefore, there exists a sequence $\{(x_k, t_k )\}$
in $(0, h(t))\times (0, \infty)$
such that $I_b(x_k,t_k)\geq \frac{\delta}{2}$ for all $k \in \mathbb{N}$, and $t_k\to \infty$ as $k\to \infty$.
Since  $0<x_k<h(t)<h_\infty<\infty$, we can choose a subsequence of $\{x_n\}$ which converges
to $x_0\in [0, h_\infty)$. Without loss of generality, we still assume $x_k\to x_0$ as $k\to \infty$.

Set $U_k(x,t)=I_b(x,t_k+t)$ and $V_k(x,t)=I_m(x,t_k+t)$ for
$x\in (0, h(t_k+t)), t\in (-t_k, \infty)$. As in \cite{FS},
 from the parabolic regularity, for $\{(U_k, V_k)\}$, we can choose a subsequence $\{(U_{k_i}, V_{k_i})\}$
 such that
$(U_{k_i}, V_{k_i})\to (\tilde U, \tilde V)$ as $i\to \infty$ and $(\tilde U, \tilde V)$ satisfies
\begin{eqnarray*} \left\{
\begin{array}{lll}
\tilde U_t-D_1 \tilde U_{xx}=\dfrac{\alpha_b\beta_b(x,t) \tilde V (N_b-\tilde U)}{N_b}-\gamma_b(x,t)
\tilde U,\; &\ 0<x<h_\infty, \ t>0,\\
\tilde V_t-D_2 \tilde V_{xx}=\dfrac{\alpha_m\beta_b(x,t) (A_m-\tilde V)\tilde U}{N_b} -d_m \tilde V,\; & 0<x<h_\infty,\ t>0.
\end{array} \right.
\end{eqnarray*}
 Recalling that $\tilde U(x_0, 0)\geq \delta/2$, therefore we derive $\tilde U>0$ in $ (0, h_\infty)\times(0, \infty)$.

 Using the similar method in proving Hopf lemma at the point $(h_\infty, 0)$ yields
that $\tilde U_x(h_\infty, 0 )\leq -\sigma_0$ for some $\sigma_0>0$.

In the meantime, since $h(t)$ is monotone increasing and bounded, for any $0<\alpha <1$ and any $\tau\geq 1$, and combining standard $L^p$ theory and the Sobolev imbedding theorem (\cite{LSU}), we
can deduce that
\begin{eqnarray}
\|I_b\|_{C^{1+\alpha,
(1+\alpha)/2}([0, h(t)]\times[\tau, \tau+1])}\leq \tilde C,\label{Bg1}
\end{eqnarray}
where the constant $\tilde C$
depends on $\alpha$, $h_0$, $\|I_{b,0}\|_{C^{2}([0, h_0])}$, $\|I_{m,0}\|_{C^{2}([0, h_0])}$ and $ h_\infty$. Note that $\tilde C$ is independent of $\tau$, by applying the free boundary conditions in (\ref{a3}), we obtain
\begin{eqnarray}
\|I_b(\cdot, t)\|_{C^{1}([0, h(t)])}\leq \hat C,\ t\geq 1 ,\label{est-2}\\
||h'||_{C^{\alpha/2}([1, +\infty))}\leq \hat C.
\label{est-1}
\end{eqnarray}

Now, since $\|h'\|_{C^{\alpha/2}([1,\infty))}\leq
\hat C$ and $h(t)$ is bounded,  we conclude that $h'(t)\to 0$ as $t\to \infty$, in the sense that
$\frac {\partial I_b}{\partial x}(h(t_k),t_k)\to 0$ as $t_k\to \infty$ by the free boundary condition. Moreover,  in view of (\ref{est-2}) we obtain
$$
\frac {\partial I_b}{\partial x}(h(t_k),t_k+0)=(U_k)_x(h(t_k),0)\to \tilde U_x(h_\infty,0),  \ \ \ \ \ \mbox{as} \ \ k\to \infty, $$
which leads to a contradiction to the fact  $\tilde U_x(h_\infty,0)\leq -\sigma_0<0$.
 Thus we have
 $$\displaystyle \lim_{t\to +\infty} \ ||I_b(\cdot,t)||_{C([0,h(t)])}=0.$$
Moreover, the above limitation indicates that for any $\varepsilon>0$, there exists a constant $T_\varepsilon>0$ such that $0\leq I_b(x,t)\leq \varepsilon$
 for $x\in [0, h(t)]$ and $t\geq T_\varepsilon$. Noting that $I_m$ satisfies
 $$\frac{\partial I_m}{\partial t}-D_2 \frac{\partial^2 I_m}{\partial x^2}\leq \frac{\alpha_m\beta_b(x,t) A_m}{N_b}\varepsilon -d_m I_m,\; 0<x<h(t),\ t\geq T_\varepsilon.$$
 Therefore $\displaystyle \limsup_{t\to +\infty} \ ||I_m(\cdot,t)||_{C([0,h(t)])}\leq \frac {\alpha_m \beta^*_b A_m}{N_b d_m}\varepsilon$, where $\beta^*_b=\sup_{[0, \infty)\times [0, T]}\beta_{b}$.
In view of the arbitrariness of $\varepsilon$,  we deduce that
$\displaystyle \lim_{t\to +\infty} \ ||I_m(\cdot,t)||_{C([0,h(t)])}=0$.
\epf

\begin{thm} Assume that $\beta_b(x,t)=\beta_b^*$ and $\gamma_b (x,t)=\gamma_b^*$. If $R_0^D([0, +\infty))\leq 1$, then $h_\infty<\infty$ and
$$\lim\limits_{t\to +\infty} ||I_b (\cdot, t)||_{C([0,h(t)])}=\lim\limits_{t\to +\infty} ||I_m (\cdot, t)||_{C([0,h(t)])}=0.$$
\label{vanish}
\end{thm}
\bpf
In this case, it is easy to check that
$$R_0^D([0, +\infty))=\ \sqrt{\frac {A_m\alpha_b\alpha_m (\beta^*_b)^2}{N_b\gamma_b^*d_m}}:=R_0.$$
We will use the energy equality to prove that $h_\infty<+\infty$. Let $k=\frac {\gamma_b^*N_b}{A_m\alpha_m \beta^*_b}$,
direct computations yield
\begin{eqnarray*}& &\frac{\textrm{d}}{\textrm{d} t}\int_{0}^{h(t)}\Big[I_b (x, t)+k I_m (x, t)\Big]\textrm{d}x\\
&=&\int_{0}^{h(t)}\Big[\frac{\partial I_b}{\partial t}+k\frac{\partial I_m}{\partial t}\Big](x, t)\textrm{d}x+h'(t)\Big[I_b +k I_m \Big](h(t), t)\\
&=&\int_{0}^{h(t)}(D_1 \frac{\partial^2 I_b}{\partial x^2}+kD_2 \frac{\partial^2 I_m}{\partial x^2})\textrm{d}x+\int_{0}^{h(t)}\, \Big[-\gamma^*_{b} I_b +\alpha_{b}\beta_{b}^*\frac {N_b-I_b}{N_{b}} I_m \Big]\textrm{d}x\\[1mm]
& &+\int_{0}^{h(t)}\, k\Big[-d_m I_m +\alpha_{m}\beta_{b}^*\frac {A_m-I_m}{N_{b}} I_b \Big]\textrm{d}x\\[1mm]
&\leq&-\frac{D_1}{\mu}h'(t)+\int_{0}^{h(t)}\Big[(-\gamma^*_{b}+k\alpha_{m}\beta_{b}^*\frac {A_m}{N_{b}})I_b +(\alpha_b\beta^*_b-k d_m) I_m\Big] \textrm{d}x\\
&=&-\frac{D_1}{\mu}h'(t)+\int_{0}^{h(t)} \alpha_b\beta^*_b(1-\frac 1{R^2_0}) I_m \textrm{d}x.
\end{eqnarray*}
Integrating from $0$ to $t\,(>0)$ gives
\begin{eqnarray}
& &\int_{0}^{h(t)}\Big[I_b+k I_m\Big](x, t)\textrm{d}x \nonumber\\
&\leq& \int ^{h(0)}_{0}\Big[I_b+k I_m\Big](x, 0)\textrm{d}x+\frac {D_1}{\mu}(h(0)-h(t))\nonumber \\
& &+\int_{0}^t\int_{0}^{h(s)}\alpha_b\beta^*_b(1-\frac 1{R^2_0}) I_m(x, s)dxds,
\quad t\geq 0.\label{k1}
\end{eqnarray}
It follows from $R_0\leq 1$ that
$$\frac {D_1}{\mu}h(t) \leq \frac {D_1}{\mu}h(0)+\int ^{h(0)}_{0}\big [I_b +k I_m\big](x, 0)\textrm{d}x$$
for $t\geq 0$, which implies that $h_\infty<\infty$. Furthermore, the vanishing of the virus
 follows easily from Lemma \ref{vanishlm}.
\epf

\begin{thm} Suppose $R_0^F(0)(:=R_0^D([0, h_0)))<1$.  Then $h_\infty<\infty$ and
$$\lim_{t\to
+\infty} \ ||I_b(\cdot, t)||_{C([0, h(t)])}=\lim_{t\to
+\infty} \ ||I_m(\cdot, t)||_{C([0, h(t)])}=0$$
provided that $||I_{m,0}(x)||_{C([0, h_0])}$
and $||I_{b,0}(x)||_{C([0, h_0])}$ are sufficiently small.
\label{vanish2}
\end{thm}
\bpf We are going to construct a suitable upper solution for problem \eqref{a3}.
Since $R_0^D([0, h_0))<1$, it follows from Theorem \ref{r0}
that there exist $\lambda_0>0$, and $\phi(x,t)>0$, $\psi(x,t)>0$ in $[0, h_0)\times[0,T]$ such that
\begin{eqnarray}
\left\{
\begin{array}{lll}
\phi_t-D_1 \phi_{xx}=\alpha_b\beta_b(x,t)\psi -\gamma_b(x,t)\phi+\lambda_0 \phi,\; & 0<x<h_0,0\leq t\leq T,\\
\psi_t-D_2 \psi_{xx}= \frac{\alpha_m\beta_b(x,t)A_m}{N_b}\phi-d_m\psi+\lambda_0 \psi,\; &0<x<h_0,0\leq t\leq T,  \\
\phi_x(0,t)=\psi_x(0,t)=0, \; & 0\leq t\leq T,\\
\phi(h_0,t)=\psi(h_0,t)=0, \; &  0\leq t\leq T,\\
\phi(x,0)=\phi(x,T), \psi(x,0)=\psi(x,T), & 0\leq x\leq h_0.
\end{array} \right.
\label{B1f1}
\end{eqnarray}

Recalling that $\phi_x(h_0,t)<0, \psi_x(h_0,t)<0$ for $t\in [0, T]$ and $\phi(x,t)>0, \psi(x,t)>0$ in $[0, h_0)\times [0, T]$. By the
regularity of $\psi$ and $\phi$, there exist  constants $C>0$ and $L>0$ such that
$$1/L\, \psi(x,t)\leq \phi(x,t) \leq L\, \psi(x,t), \; \forall (x,t)\in[0,h_0]\times[0,T].$$
$$x \phi_x(x,t)\leq C \phi(x,t),\ x \psi_x(x,t)\leq C \psi(x,t), \; \forall (x,t)\in[0,h_0]\times[0,T].$$

As in \cite{DL},  we set
$$g (t)=(1+2\delta-\delta e^{-\sigma t}),\; \xi(t)=\int_0^t g^{-2}(\tau)d\tau,\  t\geq 0,$$
$$\overline{h}(t)=h_0g(t),\; y=\frac{x}{g(t)},$$
and
$$\overline I_b(x, t)=\varepsilon e^{-\sigma t}\phi(y,\xi(t)), \ 0\leq
x\leq \overline{h}(t),\ t\geq 0,$$
$$\overline I_m(x, t)=\varepsilon e^{-\sigma t}\psi(y,\xi(t)), \ 0\leq
x\leq \overline{h}(t),\ t\geq 0,$$
where $0<\delta, \sigma, \varepsilon \ll1$ be constants, which will be chosen later.

Firstly, for any given $0<\rho\leq 1$, since
$\beta_b(x,t)$ and $\gamma_b(x,t)$ are uniformly continuous in $[0, 3h_0]\times [0,T]$ and T-periodic in $t$, then there exists
$0<\delta_0(\rho)\ll 1$ such that, for all $0<\delta\leq\delta_0(\rho)$ and $0<\sigma<1$,
we deduce that
$$ |g^{-2}(t)\beta_b(y,\xi(t))-\beta_b(x,t)|\leq\rho, \; \forall t>0,\;
 0\leq x\leq \overline{h}(t),$$
 and
 $$ |g^{-2}(t)\gamma_b(y,\xi(t))-\gamma_b(x,t)|\leq\rho, \; \forall t>0,\;
 0\leq x\leq \overline{h}(t).$$

Then, straightforward calculations yields
\begin{eqnarray*}
& &\dfrac{\partial \overline I_b}{\partial t}-D_1 \dfrac{\partial^2\overline I_b}{\partial x^2}+
\gamma_{b}(x,t)\overline I_b -\alpha_{b} \beta_{b}(x,t)\frac{(N_b-\overline I_b)}{N_b} \overline I_m \\
&\geq& -\sigma  \overline I_b+\alpha_{b}\overline I_m [\frac{1}{g^2(t)}\beta_{b}(y,\xi(t))-\beta_{b}(x,t)]\\
& &+ \overline I_b[\gamma_{b}(x,t)-\frac{1}{g^2(t)}\gamma_{b}(y,\xi(t))]+
     \lambda_0\frac{1}{g^2(t)}\overline I_b\\
& & -\varepsilon e^{-\sigma t}\frac{1}{g^2(t)}x\phi_y(y,\xi(t))\delta\sigma e^{-\sigma t}\\
&\geq & \overline I_b(-\sigma+\lambda_0\frac 1{(1+2\delta)^2}-\alpha_b \rho L-\rho-C\delta\sigma)>0,
\end{eqnarray*}
provided $0<\sigma,\rho, \delta\ll 1$ for all $t>0$ and $0\leq x\leq \overline{h}(t)$.

\begin{eqnarray*}
& &\dfrac{\partial \overline I_m}{\partial t}-D_2 \dfrac{\partial^2\overline I_b}{\partial x^2}+
d_m\overline I_m -\alpha_{m} \beta_{b}(x,t)\frac{(A_m-\overline I_m)}{N_b} \overline I_b \\
&\geq& -\sigma  \overline I_m+\alpha_{m}\overline I_b [\frac{1}{g^2(t)}\beta_{b}(y,\xi(t))-\beta_{b}(x,t)]\\
& &+ d_m\overline I_m[1-\frac{1}{g^2(t)}]+
     \lambda_0\frac{1}{g^2(t)}\overline I_m\\
& & -\varepsilon e^{-\sigma t}\frac{1}{g^2(t)}x\psi_y(y,\xi(t))\delta\sigma e^{-\sigma t}\\
&\geq & \overline I_m(-\sigma+\lambda_0\frac 1{(1+2\delta)^2}-\alpha_m \frac{A_m}{N_b}\rho L-C\delta\sigma)>0,
\end{eqnarray*}
provided $0<\sigma,\rho, \delta\ll 1$ for all $t>0$ and $0\leq x\leq \overline{h}(t)$.

Evidently, we have
$$\overline I_b(\overline{h}(t),t)=\varepsilon e^{-\sigma t}\phi(h_0, \xi(t))=0,$$
$$\overline I_m(\overline{h}(t),t)=\varepsilon e^{-\sigma t}\psi(h_0,\xi(t))=0.$$
and
$$
\begin{array}{rcl}
&&\overline{h}'(t)=\displaystyle  h_0\delta\sigma e^{-\sigma t},  \\
&& -\dfrac{\partial \overline I_b}{\partial x}(\overline{h} (t),t)=\displaystyle
   -\varepsilon \frac {1}{g(t)}\phi_y(h_0,\xi(t))e^{-\sigma t}.\\
\end{array}
   $$
If $$\varepsilon \leq \frac{h_0\delta \sigma}{\mu}\,\min_{[0,T]}\frac {-1}{\phi_y(h_0, \xi(t))},$$
 we then have
$$\overline{h}'(t)\geq -\mu \frac{\partial \overline I_b}{\partial x}(\overline{h}(t), t),$$
for $t>0$. Moreover, we now have
\begin{eqnarray*}
\left\{
\begin{array}{lll}
\dfrac{\partial \overline I_b}{\partial t}-D_1 \dfrac{\partial^2\overline I_b}{\partial x^2}\geq
\gamma_{b}(x,t)\overline I_b -\alpha_{b} \beta_{b}(x,t)\dfrac{(N_b-\overline I_b)}{N_b} \overline I_m,\; & 0<x<\overline{h}(t),\, t>0, \\
\dfrac{\partial \overline I_m}{\partial t}-D_2 \dfrac{\partial^2\overline I_b}{\partial x^2}\geq
d_m\overline I_m -\alpha_{m} \beta_{b}(x,t)\dfrac{(A_m-\overline I_m)}{N_b} \overline I_b,\; &  0<x<\overline{h}(t),\, t>0,\\
\overline I_b(x,t)=\overline I_m(x, t)=0,&x= \overline{h}(t)\, \, t>0,\\
\dfrac{\partial \overline I_b}{\partial x}(x,0)\leq 0, \dfrac{\partial \overline I_m}{\partial x}(x,0)\leq 0\; & 0<x<\overline{h}(t),\\
\overline{h}(0)> h_0, \; \overline{h}'(t)\geq -\mu \frac{\partial \overline I_b}{\partial x}(\overline{h}(t), t), & t>0.
\end{array} \right.
\end{eqnarray*}
If  $||I_{b,0}||_{L^\infty}\leq \varepsilon \min_{[0, h_0]}\phi(\frac {x}{1+\delta},0)$
and $||I_{m,0}||_{L^\infty}\leq \varepsilon \min_{[0, h_0]} \psi(\frac {x}{1+\delta},0)$, then for $x\in [0, h_0]$,
$$I_{b,0}(x)\leq \varepsilon \phi(\frac {x}{1+\delta},0)\leq \overline I_b(x, 0)$$
and
$$I_{m,0}(x)\leq \varepsilon \psi(\frac {x}{1+\delta},0)\leq \overline I_m(x, 0).$$
Then applying the comparison principle we conclude that $h(t)\leq\overline{h}(t)$ for $t>0$. It
follows that $\displaystyle h_\infty\leq \lim_{t\to\infty}
\overline{h}(t)=h_0(1+2\delta)<\infty$, and
$$\displaystyle \lim_{t\to
+\infty} \ ||I_b(\cdot, t)||_{C([0, h(t)])}= \lim_{t\to
+\infty} \ ||I_m(\cdot, t)||_{C([0, h(t)])}=0$$ by Lemma 5.2.
 \epf

Using the similar method as that in Theorem \ref{vanish2}, we can construct a suitable upper solution so that West Nile virus is vanishing when the parameter $\mu$ is sufficiently small, see also Lemma 5.10 in \cite{DL}.

 \begin{thm} Suppose $R_0^F(0)(:=R_0^D([0, h_0)))<1$.  Then there exists
$ \mu^*>0$ depending on $I_{b,0}$ and $I_{m,0}$ such that
$h_\infty<\infty$ and
$$\displaystyle \lim_{t\to
+\infty} \ ||I_b(\cdot, t)||_{C([0, h(t)])}=\lim_{t\to
+\infty} \ ||I_m(\cdot, t)||_{C([0, h(t)])}=0$$ provided that $\mu\leq \mu^*$.
\end{thm}

Next, we will give some sufficient conditions for WNv spreading. We first
exhibit that West Nile virus is spreading when $R_0^F(0)\geq 1$.
\begin{thm} If $R_0^F(0)(:=R_0^D([0, h_0)))\geq 1$, then $h_\infty=\infty$ and
$$\liminf_{t\to
+\infty} \ ||I_b(\cdot, t)||_{C([0, h(t)])}>0,\, \textrm{and}\, \liminf_{t\to
+\infty} \ ||I_m(\cdot, t)||_{C([0, h(t)])}>0,$$ that is, spreading happens.
\end{thm}
\bpf {\bf Case 1:} When $R_0^F(0):=R_0^D([0, h_0))>1$.

In this case, the following  periodic-parabolic problem
\begin{eqnarray}
\left\{
\begin{array}{lll}
\phi_t-D_1 \phi_{xx}=\alpha_b\beta_b(x,t)\psi -\gamma_b(x,t)\phi+\lambda_0 \phi,\; & 0<x<h_0,0\leq t\leq T,\\
\psi_t-D_2 \psi_{xx}= \dfrac{\alpha_m\beta_b(x,t)A_m}{N_b}\phi-d_m\psi+\lambda_0 \psi,\; &0<x<h_0,0\leq t\leq T,  \\
\phi_x(0,t)=\psi_x(0,t)=0, \; & 0\leq t\leq T,\\
\phi(h_0,t)=\psi(h_0,t)=0, \; &0\leq t\leq T,\\
\phi(x,0)=\phi(x,T), \psi(x,0)=\psi(x,T), & 0\leq x\leq h_0.
\end{array} \right.
\label{spread}
\end{eqnarray}
admits a positive solution $(\phi(x,t), \psi(x,t))$ with $||\phi||_{L^\infty}+||\psi||_{L^\infty}=1$,
which is the eigenfunction pair corresponding to the principal eigenvalue $\lambda_0<0$.
In the following, we will construct a suitable lower solution to
\eqref{a3}. For this aim, we set
$$\underline {I_b}(x,t)=\delta \phi(x,t),\quad \underline {I_m}(x,t)=\delta \psi(x,t)$$
for $0\leq x\leq h_0$, $t\geq 0$, where $\delta$ is sufficiently small which will be determined later.
Direct calculations yields
\begin{eqnarray*}
& & \dfrac{\partial \underline I_b}{\partial t}-D_1 \dfrac{\partial^2\underline I_b}{\partial x^2}+\gamma_{b}(x,t)
\underline I_b-\alpha_{b} \beta_{b}(x,t)\frac{(N_b-\underline I_b)}{N_b} \underline I_m \\
&=&\delta\phi(x,t)\big[\lambda_0+\frac{\alpha_{b}\beta_{b}(x,t)\delta\psi}{N_b}\big],\\
& &\dfrac{\partial \underline I_m}{\partial t}-D_2 \dfrac{\partial^2\underline I_m}{\partial x^2}+d_{m}\underline I_m-\alpha_{m} \beta_{b}(x,t) \frac{(A_{m}- \underline I_m)}{N_b} \underline I_b\\
&=&\delta\psi(x,t)\big[\lambda_0+\frac{\alpha_{m}\beta_{b}(x,t)\delta\phi}{N_b}\big],
\end{eqnarray*}
for  $0<x<h_0$ and $0\leq t\leq T$.

Recalling $\lambda_0<0$, we can choose $\delta$ sufficiently small such that
\begin{eqnarray*}
\left\{
\begin{array}{lll}
\dfrac{\partial \underline I_b}{\partial t}-D_1 \dfrac{\partial^2\underline I_b}{\partial x^2}
\leq -\gamma_{b}(x,t)\underline I_b+\alpha_{b} \beta_{b}(x,t)\dfrac{(N_b-\underline I_b)}{N_b}
 \underline I_m, \; &0<x<h_0,\, 0\leq t\leq T, \\
\dfrac{\partial \underline I_m}{\partial t}-D_2 \dfrac{\partial^2\underline I_m}{\partial x^2}\leq -d_{m}\underline I_m+\alpha_{m} \beta_{b}(x,t) \dfrac{(A_{m}- \underline I_m)}{N_b} \underline I_b,
\; &0<x<h_0,\, 0\leq t\leq T, \\
\underline I_b(x,t)=\underline I_m(x, t)=0, \; &x= h_0,\,0\leq t\leq T, \\
0=h'_0\leq -\mu\frac{\partial \underline I_b}{\partial x}(h_0, t), \; &t>0,\\
\underline {I_b}(x,0)\leq I_{b,0}(x),\ \underline{I_m}(x,0)\leq I_{m,0}(x),\; &0\leq x\leq h_0.
\end{array} \right.
\end{eqnarray*}
Therefore, by applying comparison principle, we derive that $I_b(x,t)\geq\underline I_b(x,t)$ and  $I_m(x,t)\geq\underline I_m(x,t)$
in $[0, h_0]\times [0,T]$. It follows that $\displaystyle \liminf_{t\to
+\infty} \|I_b(\cdot, t)||_{C([0, h(t)])}\geq \delta \phi(0)>0$ and $\displaystyle \liminf_{t\to
+\infty} \|I_m(\cdot, t)||_{C([0, h(t)])}\geq \delta \psi(0)>0$, therefore $h_\infty=+\infty$ by Lemma \ref{vanishlm}.

{\bf Case 2:} When $R_0^F(0):=R_0^D([0, h_0))=1$.

If $R_0^F(0):=R_0^D([0, h_0))=1$, then for any positive time $t_0$, we deduce that $h(t_0)>h(0)=h_0$ by Theorem \ref{exist}, therefore $R_0^D([0, h(t_0)))>R_0^D([0, h_0))=1$ by the monotonicity in Lemma \ref{proper}.
Substituting the initial time $0$ by the positive time $t_0$, we derive $h_\infty=+\infty$ as Case 1.
 \epf

\begin{rmk} From the above proof, one can see that the spreading scenario will happens if there exists $t_0\geq 0$ such that $R_0^F(t_0)\geq 1$. Moreover, if $R^D_0([0, +\infty))>1$, the condition is sufficient and necessary. In fact, if $R_0^F(t)< 1$ for any $t\geq 0$, we then have $h_\infty<+\infty$ and vanishing happens.
\end{rmk}

In the following, we explore the long time asymptotic behavior of the solution to problem \eqref{a3} when the spreading happens.
\begin{thm} \label{asymp} If $R_0^F(t_0)(:=R_0^D([0, h(t_0)))\geq 1$ for some $t_0\geq 0$, then $h_\infty=+\infty$ and
\begin{eqnarray}
\left\{
\begin{array}{ll}
\underline{U}(x,t)\leq\liminf\limits_{n\longrightarrow\infty}I_b(x,t+nT)\leq
\limsup\limits_{n\longrightarrow\infty}I_b(x,t+nT)\leq\overline{U}(x,t),\\
\underline{V}(x,t)\leq\liminf\limits_{n\longrightarrow\infty}I_m(x,t+nT)\leq
\limsup\limits_{n\longrightarrow\infty}I_m(x,t+nT)\leq\overline{V}(x,t),
\end{array} \right.
\label{bd5}
\end{eqnarray}
uniformly holds in any compact subset of $[0, \infty)\times[0,T]$, where $(\overline{U}, \overline{V})$ and
$(\underline{U},\underline{V})$ are the maximal and the minimal positive periodic solutions of
the corresponding  T-periodic boundary value
problem \eqref{b1} in half space.
\end{thm}
\bpf It follows from Theorem \ref{theb5} that, for any positive constant $l>h(t_0)$, $(\overline{U}_l, \overline {V}_l)$ converges to the $(\overline{U}, \overline{V})$ locally uniformly in $[0,+\infty)\times[0,T]$ as $l\rightarrow\infty$, which is the maximal positive periodic solutions of problem \eqref{b1}.
Therefore, for any given $L_1\in (0, +\infty)$, $(\overline{U}_l, \overline {V}_l)$ converges to $(\overline{U}, \overline{V})$ uniformly on $[0, L_1]\times[0,T]$.
Hence, for any $\varepsilon>0$, there exists a positive constant $l_\varepsilon> \max\{L_1, h(t_0)\}$
such that
\begin{equation}\label{ineq}
(\overline{U}_{l_\varepsilon}, \overline {V}_{l_\varepsilon})\leq (\overline{U}+\varepsilon, \overline{V}+\varepsilon)
\end{equation}
 on $[0, L_1]\times [0,T]$.

On the other hand, for the above $l_\varepsilon$, there exists a constant $T_\varepsilon$ such that
$$R_0^F(T_\varepsilon)\geq R_0^D([0,l_\varepsilon))> R_0^D([0,h(t_0))>1,$$
  and comparing $(N_b, A_m)$ with $(I_b, I_m)$ in $[0, l_\varepsilon]\times [T_\varepsilon, \infty)$,
  it is easy to see that $(N_b, A_m)$  and $(I_b, I_m)$ are ordered upper and lower solution of the system in $[0, l_\varepsilon]\times [T_\varepsilon, \infty)$. By the same procedure as in the proof of Theorem \ref{theb5}, we can deduce that
$$ \limsup_{n\rightarrow\infty}(I_b, I_m)(x,t+nT)\leq (\overline{U}_{l_\varepsilon}, \overline {V}_{l_\varepsilon})(x,t)$$
in $[0, l_\varepsilon]\times[0,T]$, which together with \eqref{ineq} implies
$$ \limsup_{n\rightarrow\infty}(I_b, I_m)(x,t+nT)\leq (\overline{U}+\varepsilon, \overline{V}+\varepsilon)$$
on $[0, L_1]\times [0,T]$. Then we obtain
$$ \limsup_{n\rightarrow\infty}(I_b, I_m)(x,t+nT)\leq (\overline{U}, \overline{V})$$
uniformly in $[0, L_1]\times [0, T]$ due to the arbitrariness of $\varepsilon$.
The remaining two inequalities can be proved similarly.
\epf

\end{document}